\documentclass[pdflatex,sn-mathphys-num]{sn-jnl}

\usepackage{graphicx}%
\usepackage{multirow}%
\usepackage{amsmath,amssymb,amsfonts}%
\usepackage{amsthm}%
\usepackage{mathrsfs}%
\usepackage[title]{appendix}%
\usepackage{xcolor}%
\usepackage{textcomp}%
\usepackage{manyfoot}%
\usepackage{booktabs}%
\usepackage{algorithm}%
\usepackage{algorithmicx}%
\usepackage{algpseudocode}%
\usepackage{listings}%

\theoremstyle{thmstyleone}%
\newtheorem{theorem}{Theorem}
\newtheorem{coro}[theorem]{Corollary}
\theoremstyle{thmstyletwo}%

\theoremstyle{thmstylethree}%

\begin{document}

\title[Article Title]{The Optimal Strategy for Playing \textit{Lucky 13}}

\author*[1]{\fnm{Steven} \sur{Berger}}\email{srb420@lehigh.edu}

\author[1]{\fnm{Daniel} \sur{Conus}}\email{dac311@lehigh.edu}

\affil*[1]{\orgdiv{Department of Mathematics}, \orgname{Lehigh University}, \orgaddress{\street{Chandler-Ullmann Hall, 17 Memorial Drive East}, \city{Bethlehem}, \postcode{18015}, \state{PA}, \country{United States}}}


\abstract{The game show \textit{Lucky 13} differs from other television game shows in that contestants are required to place a bet on their own knowledge of trivia by selecting a range that contains the number of questions that they answered correctly.  We present several models for this game show using binomial random variables and generate tables outlining the optimal range the player should select based on maximization of two different utility functions. After analyzing the decisions made by some actual contestants on this show, we present a numerical simulation for how many questions an average player is expected to answer correctly based on question categories observed for twelve sample contestants.   }


\maketitle

\section{Introduction}\label{sec1}
Strategic decision making in games provides a unique and entertaining application of probability theory. There are strategies available for most common casino games.  For example, in \citep{Bib3}, the optimal strategy for blackjack is discussed depending on whether or not the player wants to count cards. Many television game shows also provide scenarios in which strategic play must be applied.  Perhaps the most famous example is the Monty Hall Problem commonly seen in almost every introductory probability book.  The original version of the game and some of its variations are discussed in \citep{Bib4}.  In television game shows where contestants are asked to answer trivia questions, one can usually win more money by answering more questions correctly.  For example, in the popular game show ``Who wants to be a millionaire" contestants were asked a series of 14 or 15 questions and were allowed to use a set of lifelines on questions or walk away from a question if they were not confident. Using stochastic control techniques, Dalang and Bernyk \citep{Bib1} have created tables to show the optimal lifeline strategy, upper bounds on the expected payoff when the optimal strategy is followed, and the upper bounds on the probability of winning the million dollars given a certain number of questions have been correctly answered. However, in the 2024 ABC television game show \textit{Lucky 13}, contestants not only had to answer the questions correctly, but also correctly predict how many correct answers they gave.  In this scenario, your ability to evaluate your own knowledge is just as important as providing the maximum number of correct responses.  \\
\\

\section{Game Description}
In \textit{Lucky 13}, each contestant answers 13 trivia questions where every answer is either ``True'' or ``False".  After answering all 13 questions but before the correct answers are revealed, the contestants are asked to select a ``Lucky Range" from Table \ref{Tab1} containing the number of questions out of 13 they answered correctly.  A contestant wins the prize amount from this table if they correctly predict their Lucky Range before the correct answers are revealed.  Contestants are told they will win nothing if their actual number of correct responses is under or over their predicted Lucky Range.  Note that a contestant cannot predict 0 questions were answered correctly.  Additionally, contestants who did not pick Lucky Range 13 are given the opportunity to predict their exact number of correct responses (called the ``Lucky Number") from their Lucky Range for an additional \$25,000.  For example, a contestant who predicts Lucky Range 7-9 and Lucky Number 9 would win \$25,000 if he got 8 questions correct, \$50,000 if he got 9 correct, and \$0 if he got 10 correct.  Once the contestant has selected a Lucky Range and Lucky Number, correct answers are revealed one by one in an order chosen by the show's producer.  At some point before all correct answers have been revealed, the contestant is offered guaranteed cash to immediately stop playing.  The contestant may accept this offer or continue playing.  \\

\begin{table}[h]
\caption{Table of Lucky Ranges and Cash Prizes}\label{Tab1}%
\begin{tabular}{@{}ll@{}}
\toprule
Lucky Range & Prize \\
\midrule
1-3    & \$5,000   \\
4-6    & \$15,000   \\
7-9    & \$25,000  \\
10-12 & \$100,000\\
13 & \$1,000,000 \\
\botrule
\end{tabular}
\end{table}

Before playing this game, a contestant probably has already decided how much value each cash prize has.  For example, a wealthy businessman who regularly plays trivia games would probably place a higher personal value on winning \$1,000,000 than on winning any lower prize.  However, other contestants may simply be interested in winning some amount of money and would be less willing to engage in risky play.  Mathematically, this personal value of the cash prizes can be quantified by a utility function that maps each prize to a real number where higher values correspond to a more desirable personal value.  Appropriate mathematical axioms needed to define a utility function on a discrete set are given in \citep{Bib2}.  In this paper, we apply our models to two different utility functions.  Our models can be applied to other utility functions should a player choose a different one than the ones discussed here.  

\section{Model with Two Question Categories}
In order to determine a strategy for correctly selecting the Lucky Range and Lucky Number, we consider a model where each of the 13 questions can be categorized into one of two different types.  A question is ``Sure" (S) if the contestant gives the correct answer with probability 1.  A question is ``Guess" (G) if the contestant will give the correct answer with probability 0.5. Let $N_S$ be the number of Sure questions (out of 13) a contestant has.  Let $N_G$ be the number of Guess questions that the contestant answers correctly.  These are the extreme cases since we assume either a contestant knows the answer for sure or randomly guesses.  We also assume that the contestant is able to decide which category is appropriate for each question.  In the next section we introduce a model that takes partial knowledge of a question into account.  Clearly, $N_S$ is not random but $N_G$ is a random variable that follows a binomial distribution.  More specifically, $N_G \sim \text{Binom}(13-N_S, 0.5)$.  Additionally, define the random variable $N$ to be the total number of correct answers out of 13, so that $N = N_S+N_G$.  \\
\\
Now we will examine several cases.  First, assume we have Contestant $X$ playing this game who guessed on each of the 13 questions (so $N_S = 0 $).  For Contestant $X$, a basic computation shows that $\text{E}[N] = \text{E}[N_G] = \frac{13}{2}=6.5$ and $\text{Var}[N] = \text{Var}[N_G] =\frac{13}{4} = 3.25$.  We can also easily calculate the probability for $N$ to fall in each Lucky Range:  
$$\text{P}(1 \leq N \leq 3) = \left[ {13 \choose 1} + {13 \choose 2}+ {13 \choose 3} \right]\cdot \left(\frac{1}{2} \right)^{13} \approx 0.0460$$
$$\text{P}(4 \leq N \leq 6) = \left[ {13 \choose 4} + {13 \choose 5}+ {13 \choose 6} \right]\cdot \left(\frac{1}{2} \right)^{13} \approx 0.4539$$
$$\text{P}(7 \leq N \leq 9) = \left[ {13 \choose 7} + {13 \choose 8}+ {13 \choose 9} \right]\cdot \left(\frac{1}{2} \right)^{13} \approx 0.4539$$
$$\text{P}(10 \leq N \leq 12) = \left[ {13 \choose 10} + {13 \choose 11}+ {13 \choose 12} \right]\cdot \left(\frac{1}{2} \right)^{13} \approx 0.0460$$
$$\text{P}(N = 13) = \left(\frac{1}{2} \right)^{13} \approx 0.0001$$
Table \ref{Tab2} gives the expected winnings for each Lucky Range.  
\begin{table}[h]
\caption{Table of Lucky Ranges and Expected Winnings when $N_S = 0$}\label{Tab2}%
\begin{tabular}{@{}ll@{}}
\toprule
Lucky Range & Expected Winnings \\
\midrule
1-3    & \$5,000(0.0460) = \$230   \\
4-6    & \$15,000(0.4539) = \$6,808.5   \\
7-9    & \$25,000(0.4539) = \$11,347.5  \\
10-12 & \$100,000(0.0460) = \$4,600\\
13 & \$1,000,000(0.0001) = \$100 \\
\botrule
\end{tabular}
\end{table}
It is clear from the probability distribution of $N$ that Contestant $X$ should choose either the Lucky Range 4-6 or 7-9.  Since the mathematical expectation for his number of correct responses is 6.5 and the probability of falling in Lucky Range 4-6 is equal to the probability of falling in Lucky Range 7-9, he should choose 7-9 as that has a higher cash prize than 4-6.  However, if Contestant $X$ strongly feels like he only answered 6 or less correctly, it may be advantageous for him to select Lucky Range 4-6.  \\
\\
Now suppose that Contestant $Y$ is playing and has 3 Sure questions.  In this case, $\text{E}[N] = 3+ \text{E}[N_G] = 3+\frac{10}{2}= 8 $ and $\text{Var}[N] = \frac{10}{4} = 2.5$.  The probability distribution for $N$ is now:
$$\text{P}(1 \leq N \leq 3) =  {10 \choose 0}\cdot \left(\frac{1}{2} \right)^{10} \approx 0.0009$$
$$\text{P}(4 \leq N \leq 6) = \left[ {10 \choose 1} + {10 \choose 2}+ {10 \choose 3} \right]\cdot \left(\frac{1}{2} \right)^{10} \approx 0.1709$$
$$\text{P}(7 \leq N \leq 9) = \left[ {10 \choose 4} + {10 \choose 5}+ {10 \choose 6} \right]\cdot \left(\frac{1}{2} \right)^{10} \approx 0.6563$$
$$\text{P}(10 \leq N \leq 12) = \left[ {10 \choose 7} + {10 \choose 8}+ {10 \choose 9} \right]\cdot \left(\frac{1}{2} \right)^{10} \approx 0.1709$$
$$\text{P}(N = 13) = {10 \choose 10}\left(\frac{1}{2} \right)^{10} \approx 0.0009$$
Table \ref{Tab3} gives Contestant $Y$'s expected winnings for each Lucky Range. 
\begin{table}[h]
\caption{Table of Lucky Ranges and Expected Winnings when $N_S = 3$}\label{Tab3}%
\begin{tabular}{@{}ll@{}}
\toprule
Lucky Range & Expected Winnings \\
\midrule
1-3    & \$5,000(0.0009) = \$4.5   \\
4-6    & \$15,000(0.1709) = \$2,563.5   \\
7-9    & \$25,000(0.6563) = \$16,407.5  \\
10-12 & \$100,000(0.1709) = \$17,090\\
13 & \$1,000,000(0.0009) = \$900 \\
\botrule
\end{tabular}
\end{table}
In this case, Contestant $Y$ will need to apply some personal judgment in order to select his Lucky Range.  If he wants to maximize the probability he wins some money, he should choose Lucky Range 7-9 because this one has the greatest probability.  However, we see from Table \ref{Tab3} that the mathematical expectation for his winnings is higher if he picks Lucky Range 10-12.  Although it is not possible for Contestant $Y$ to win exactly \$17,090 (the mathematical expected winning), selecting Lucky Range 10-12 would be the better strategy in the long term if Contestant $Y$ were playing this game multiple times and always was sure about 3 answers.  However, for the single time he must make his decision, he needs to evaluate how much the cash prizes are worth to him.  If Contestant $Y$ only wants to win a large sum of money and likes to take risks, he should select 10-12.  On the other hand, if he is a more risk-averse player, he should select 7-9.  This example shows that for this game show, maximizing the probability of winning does not necessarily lead to selecting the same Lucky Range as a player who chooses to maximize his potential cash winnings.

\begin{table}[h]
\caption{Strategy for Selecting Lucky Range and Number for 2 Question Categories}\label{Tab4}
\begin{tabular*}{\textwidth}{@{\extracolsep\fill}lcccc}
\toprule%
& \multicolumn{2}{@{}c@{}}{Maximizing Win Probability} & \multicolumn{2}{@{}c@{}}{Maximizing Winnings} \\ \cmidrule{2-3} \cmidrule{4-5}
Sure & Lucky Range & Lucky Number & Lucky Range & Lucky Number  \\
\midrule
0  & 7-9\footnotemark[1]  & 7 & 7-9 & 7 \\
1  & 7-9 & 7  & 7-9  & 7\\
2  & 7-9 & 8  & 7-9  & 8\\
3  & 7-9 & 8  & 10-12  & 10\\
4  & 7-9 & 9  & 10-12  & 10\\
5  & 7-9 & 9  & 10-12  & 10\\
6  & 10-12\footnotemark[2] & 10  & 10-12  & 10\\
7  & 10-12 & 10  & 10-12  & 10\\
8  & 10-12 & 11  & 10-12  & 11\\
9  & 10-12 & 11  & 10-12  & 11\\
10  & 10-12 & 12  & 13  & NA\\
11  & 10-12 & 12  & 13  & NA\\
12  & 13\footnotemark[3] & NA  & 13  & NA\\
13  & 13 & NA  & 13  & NA\\

\botrule

\end{tabular*}
\footnotetext[1]{Lucky Range 4-6 with Lucky Number 6 also has the same probability.  }
\footnotetext[2]{Lucky Range 7-9 with Lucky Number 9 also has the same probability.  }
\footnotetext[3]{Lucky Range 10-12 with Lucky Number 12 also has the same probability.  }
\end{table}
Table \ref{Tab4} gives the appropriate Lucky Range and Lucky Number to select based on your decision to maximize either the probability that you win some money or your expected winnings.  We see that maximizing expected winnings results in more risk as contestants may play for a million dollars if they only guessed on 3 questions whereas players maximizing their win probability should select Lucky Range 10-12 in that case.  This table represents two possible extremes for risk tolerance and some players may be willing to take some level of risk that lies in between these two extremes.  Note that maximizing the probability of winning is represented by a constant utility function while maximizing expected winnings is represented by a linear utility function.  A randomly selected player may have selected a different utility function but would still be able to use Table \ref{Tab4} to guide his decision on selecting the Lucky Range and Number.

\section{Model with Three Question Categories}
While the model with two question categories may be suitable in some scenarios, other players may have some partial knowledge about a question without knowing the answer for sure.  In addition to ``Sure" and ``Guess" questions (defined above), we now allow some questions to be in a third category ``Unsure" (U) if the contestant will give the correct answer with probability 0.75.  Let $N_S,N_U,$ and $N_G$ be the number of correct responses to Sure, Unsure, and Guess questions respectively.  As before, let $N$ denote the total number of correct responses so that $N=N_S+N_U+N_G$.  It is clear that $N_U$ and $N_G$ both follow the binomial distribution with parameter $p = 0.75$ and $0.5$ respectively.  One may obtain the same results as in the two category model by setting $N_U=0$. 
 \\
\\
We now consider the case of Contestant $Z$ who has 10 Sure questions, 2 Unsure questions, and 1 Guess question.  In this case, $\text{E}[N] =10+ \text{E}[N_S] + \text{E}[N_G] =10 + \frac{6}{4} + \frac{1}{2}  = 12$.  The probability that $N$ lands in different Lucky Ranges is as follows:
$$\text{P}(1 \leq N \leq 3) = \text{P}(4 \leq N \leq 6) = \text{P}(7 \leq N \leq 9) = 0$$
$$\text{P}(N = 13) = \text{P}(N_U = 2, N_G = 1) = \left(\frac{3}{4}\right)^2\cdot \left(\frac{1}{2}\right) = \frac{9}{32} \approx 0.2813$$
$$\text{P}(10 \leq N \leq 12) = 1 - \text{P}(N = 13) - \text{P}(N \leq 9) = 1- \frac{9}{32} - 0 = \frac{23}{32} \approx 0.7188$$

\begin{table}[h]
\caption{Table of Lucky Ranges and Expected Winnings for 10 Sure, 2 Unsure, and 1 Guess}\label{Tab5}%
\begin{tabular}{@{}ll@{}}
\toprule
Lucky Range & Expected Winnings \\
\midrule
1-3    & \$0  \\
4-6    & \$0  \\
7-9    & \$0 \\
10-12 & \$100,000(0.7188) = \$71,880\\
13 & \$1,000,000(0.2813) = \$281,300 \\
\botrule
\end{tabular}
\end{table}

Table \ref{Tab5} shows that Contestant $Z$ should play for a million dollars if he is willing to tolerate substantial risk, but otherwise should play for \$100,000 and select Lucky Range 10-12.  We can also calculate which Lucky Number Contestant $Z$ should select.  \\
\\
$$\text{P}(N=10) = \text{P}(N_U = 0, N_G = 0) = {2 \choose 0}\left(\frac{1}{4} \right)^2 \cdot {1 \choose 1}\left( \frac{1}{2}\right)^1 = \frac{1}{32}$$

$$\text{P}(N = 11) = \text{P}(N_U =0, N_G =1) + \text{P}(N_U=1, N_G=0)$$
$$= {2 \choose 0}\left(\frac{1}{4}\right)^2\cdot {1 \choose 1}\left(\frac{1}{2}\right)^1 + {2 \choose 1}\left(\frac{3}{4} \right)^1\left(\frac{1}{4} \right)^1\cdot {1 \choose 0}\left(\frac{1}{2}\right)^1 = \frac{7}{32}$$

$$\text{P}(N = 12) = \text{P}(N_U = 1, N_G = 1)+\text{P}(N_U = 2, N_G = 0)$$
$$ = {2 \choose 1}\left(\frac{3}{4} \right)^1\left(\frac{1}{4} \right)^1\cdot {1 \choose 1}\left(\frac{1}{2}\right)^1+ {2 \choose 2}\left(\frac{3}{4} \right)^2\cdot {1 \choose 0}\left( \frac{1}{2}\right)^1 = \frac{15}{32}$$
Clearly, Contestant $Z$ should select Lucky Number 12.  \\
\\

There are 105 different combinations of Sure, Unsure, and Guess questions.  We present the optimal decision based on the two different utility functions for choosing a Lucky Number and Lucky Range for selected combinations in Table \ref{Tab6}. The strategy for all combinations of $S,U, $ and $G$ can be found in Appendix A.  
 We see that it is never optimal to select Lucky Range 1-3 and risk-averse players should only select Lucky Range 4-6 if they guess on all 13 questions.  Maximizing expected winnings corresponds to more risky play when a contestant's expected number of correct answers is between one of the Lucky Ranges in which case they choose the range with higher cash prize (even if there is a lower probability of being in that range).  

\begin{table}[h]
\caption{Selected Lucky Ranges and Numbers for 3 Question Categories}\label{Tab6}
\begin{tabular*}{\textwidth}{@{\extracolsep\fill}lcccc}
\toprule%
& \multicolumn{2}{@{}c@{}}{Maximizing Win Probability} & \multicolumn{2}{@{}c@{}}{Maximizing Expected Winnings} \\ \cmidrule{2-3} \cmidrule{4-5}
S/U/G & Lucky Range & Lucky Number & Lucky Range & Lucky Number  \\
\midrule
0/0/13  & 7-9\footnotemark[1]  & 7 & 7-9 & 7 \\
0/1/12  & 7-9  & 7 & 7-9 & 7 \\
1/1/11 & 7-9 & 7 & 7-9 & 7 \\
3/8/2 & 10-12 & 10 & 10-12 & 10\\
5/6/2 & 10-12 & 11 & 10-12 &11\\
5/7/1 & 10-12 & 11 & 10-12 & 11\\
5/8/0 & 10-12 & 11 & 13 & NA\\
6/5/2 & 10-12 & 11 & 10-12 & 11\\
6/6/1 & 10-12 & 11 & 13 & NA\\
6/7/0 & 10-12 & 11 & 13 & NA\\
7/0/6 & 10-12 & 10 & 10-12 & 10\\
7/1/5 & 10-12 & 10 & 10-12 & 10\\
7/2/4 & 10-12 & 11 & 10-12 & 11\\
7/3/3 & 10-12 & 11 & 10-12 & 11\\
7/4/2 & 10-12 & 11 & 10-12 & 11\\
7/5/1 & 10-12 & 11 & 13 & NA\\
7/6/0 & 10-12 & 12 & 13 & NA\\
8/2/3 & 10-12 & 11 & 10-12 & 11\\
8/3/2 & 10-12 & 11 & 13 & NA\\
9/4/0 & 10-12 & 12 & 13 & NA\\
10/3/0 & 10-12 & 12 & 13 & NA\\
10/2/1 & 10-12 & 12 & 13 & NA\\
11/1/1 & 10-12 & 12 & 13 & NA\\
12/1/0 & 13 & NA & 13 & NA\\

\botrule
\end{tabular*}
\footnotetext[1]{Lucky Range 4-6 with Lucky Number 6 also has the same probability.  }

\end{table}

\section{Monte Carlo Simulation}
\begin{figure}[h!]
    \centering
    \includegraphics[width=0.75\linewidth]{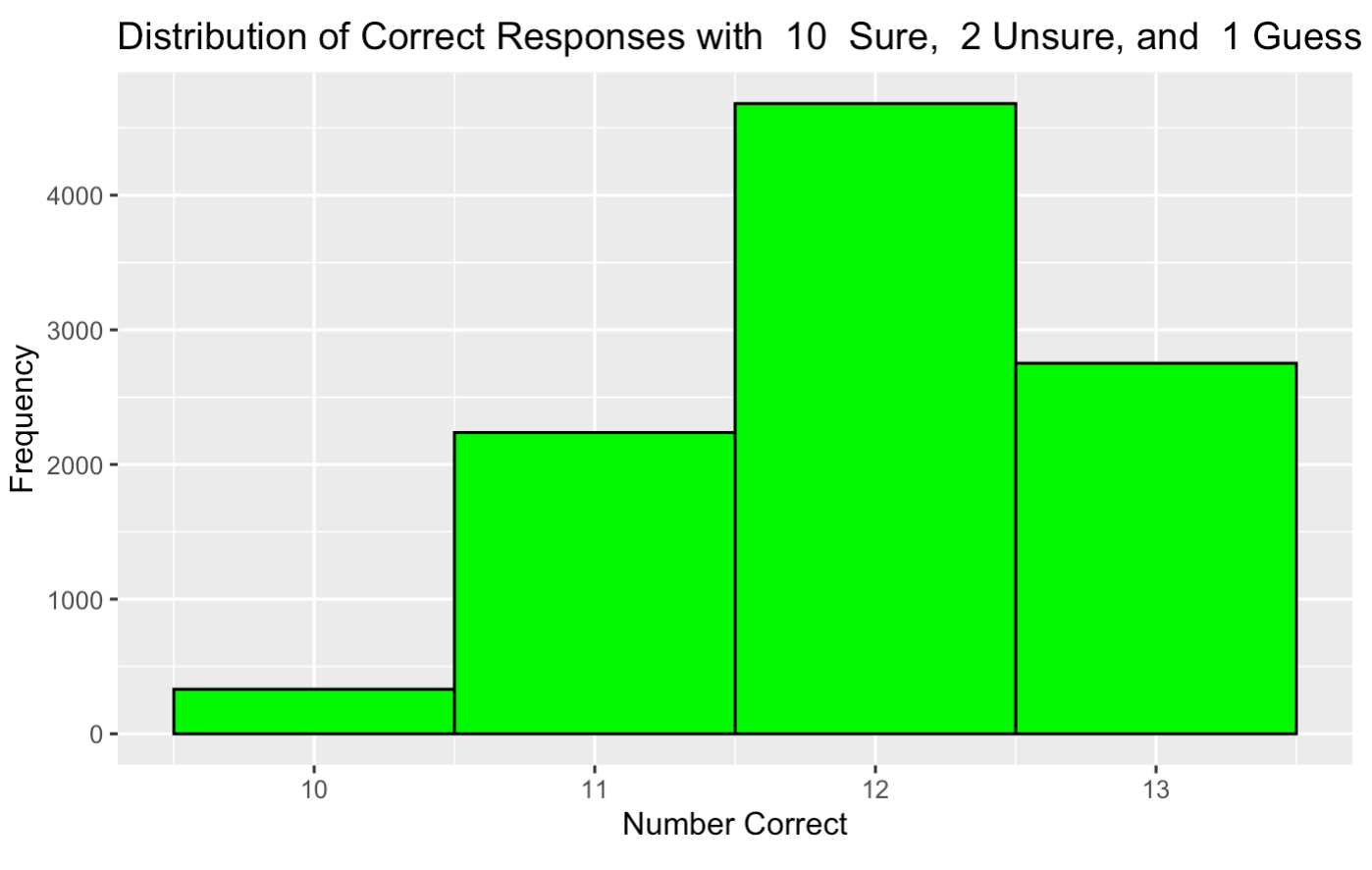}
    \caption{10,000 Monte Carlo Simulation Results for Contestant $Z$}
    \label{Fig1}
\end{figure}

The results shown in Table \ref{Tab6} (calculated using the binomial distribution) were checked by Monte Carlo simulation.  We simulate and record $N = S + \text{Binom}(U,0.75) + \text{Binom}(G,0.5)$ for 10,000 times.  In the case described above for Contestant $Z$, simulation gives the following probability distribution which agrees with the exact calculation given previously.
$$\text{P}(1 \leq N \leq 3) = \text{P}(4 \leq N \leq 6) = \text{P}(7 \leq N \leq 9) = 0$$

$$\text{P}(N = 13)   \approx 0.2844$$

$$\text{P}(10 \leq N \leq 12) \approx 0.7156$$
\\
Figure \ref{Fig1} shows the histogram containing the simulated values of $N$ for Contestant $Z$.  We see that he is most likely to get 12 questions correct which confirms the choice presented earlier of Lucky Range 10-12 with Lucky Number 12.  

\section{Case Study 1}
\begin{figure}[h!]
    \centering
    \includegraphics[width=0.8\linewidth]{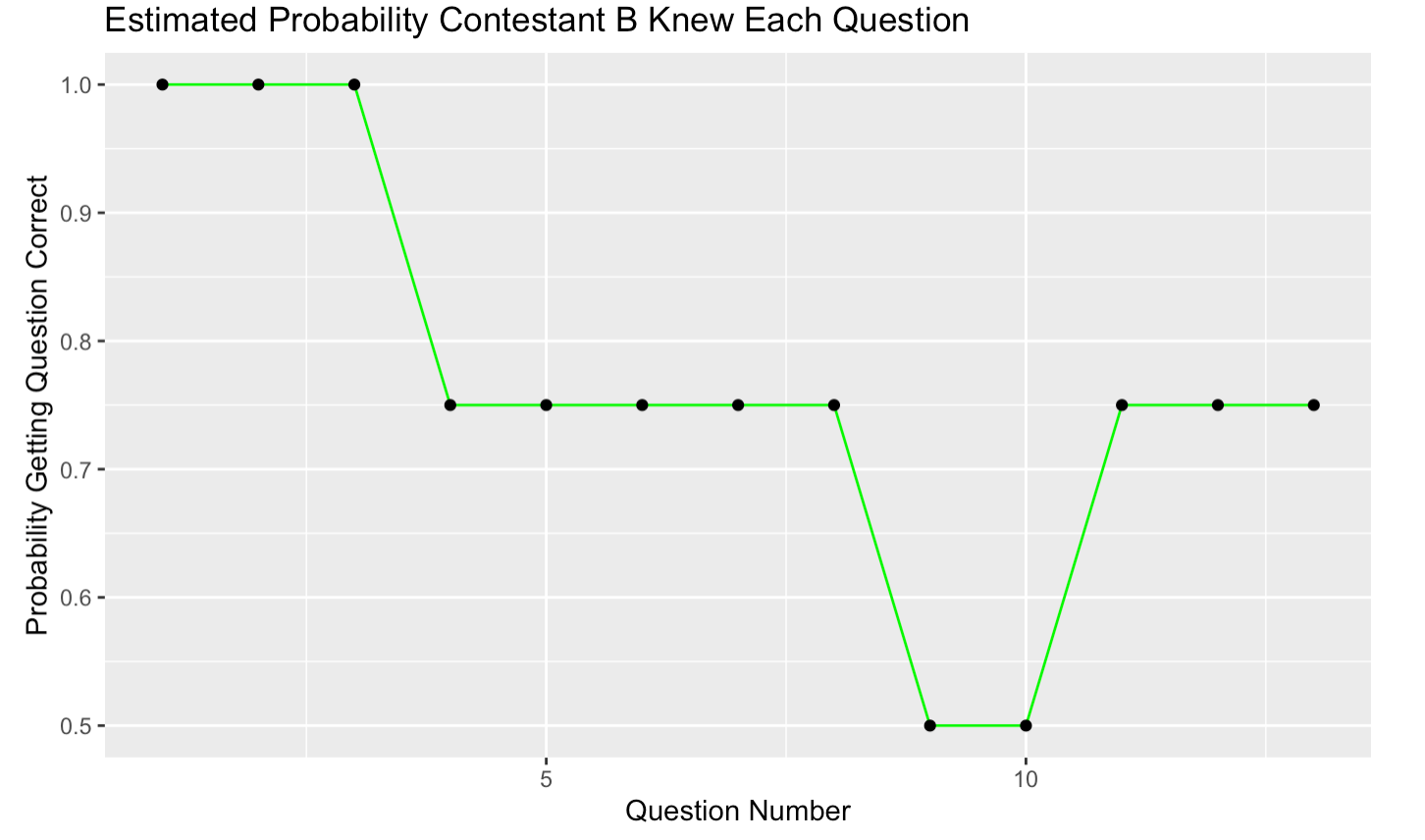}
    \caption{Estimated Probabilities of Correct Answers for Contestant $B$}
    \label{Fig2}
\end{figure}

\begin{figure}[h!]
    \centering
    \includegraphics[width=0.80\linewidth]{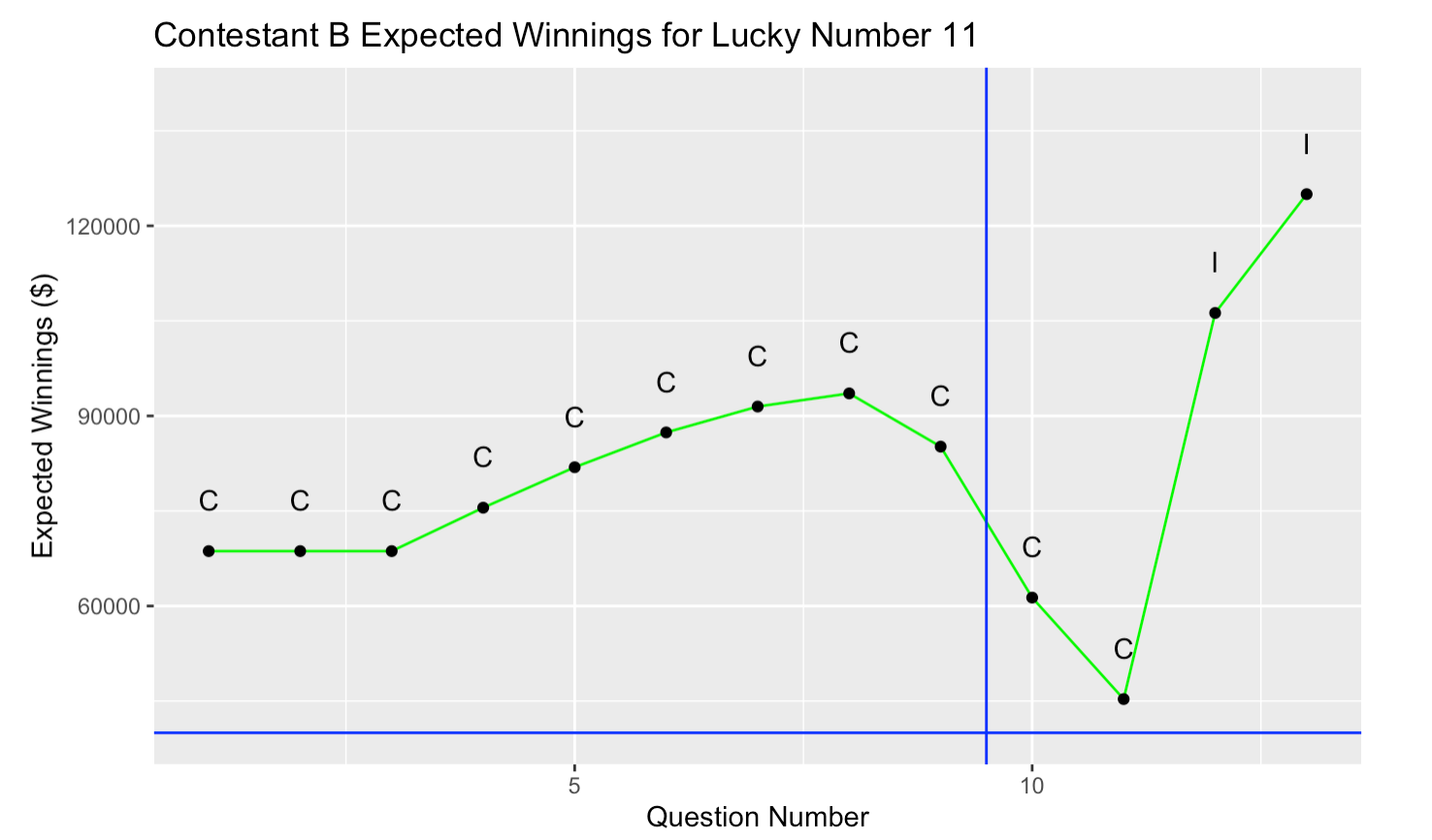}
    \caption{Expected Winnings for Contestant $B$}
    \label{Fig3}
\end{figure}

Now, we will discuss the game of an actual contestant on this show who we will refer to as Contestant $B$.  After watching Contestant $B$ answer the 13 questions, the authors determined (using our subjective judgment) that he had $S=3, U=8, $ and $G = 2$.  Figure \ref{Fig2} shows the order in which the producers chose to reveal his answers. 
 Contestant $B$ selected Lucky Range 10-12 and Lucky Number 11.  According to Table \ref{Tab6}, this was the correct decision for Lucky Range, but the incorrect decision for Lucky Number (regardless of utility function).  With $U=8$ and $G =2$ there are 27 possible combinations for $N_U$ and $N_G$.  For example, one possible combination is $\text{P}(N_U = 8, N_G = 0) \approx 0.025$ which would result in Contestant $B$ hitting his Lucky Number since $S=N_S=3$.  By summing the probabilities of relevant outcomes, we see that before any answers are revealed, Contestant $B$ has approximately a 63\% chance of hitting his Lucky Range and a 24\% chance of hitting his Lucky Number.  This means that Contestant $B$'s mathematical expectation for winnings is approximately $\$100,000(0.63)+\$25,000(0.24) = \$69,000$.  Using more decimal places, a computer calculates the mathematical expectation to be \$68,665.41.  Since the first answer revealed was a Sure question that Contestant $B$ answered correctly, this mathematical expectation does not change after one answer is revealed. Calculating this mathematical expectation after the reveal of each answer is important as it allows Contestant $B$ to evaluate his offer to quit when it is announced.  In this game, the producers of the show offered Contestant $B$, \$40,000 to stop playing after nine answers were revealed.  At this point in the game, Contestant $B$ had 9 correct answers with 3 Unsure questions and 1 Guess question remaining.  His mathematical expectation for winnings after nine answers revealed is \$85,156.25.  Naturally, Contestant $B$ rejected the offer from the producers and kept playing.  With all nine questions correct so far, the probability Contestant $B$ hits his Lucky Number is 
 $$\text{P}(N=11) = \text{P}(N_U = 2, N_G=0)+\text{P}(N_U=1,N_G=1) = \frac{9}{32}$$
 and the probability he is in his Lucky Range is
 $$\text{P}(\text{in Lucky Range }) = \text{P}(N = 10) + \text{P}(N=11) + \text{P}(N = 12)$$
 $$=\text{P}(N_U = 1, N_G = 0)+\text{P}(N_U = 0, N_G =1)+ \frac{9}{32} +$$
 $$\text{P}(N_U = 3, N_G = 0)+\text{P}(N_U = 2, N_G = 1)$$
 $$ = \frac{25}{32}$$
 At the end of the game, Contestant $B$ had only two incorrect responses (both given for Unsure questions) so he won the prizes for hitting both his Lucky Range and Lucky Number.  Figure $\ref{Fig3}$ shows the trajectory of Contestant $B$'s expected winnings after each correct answer is revealed.  The solid parallel lines intersecting at the point $(9.5,40,000)$ show the point at which his offer to quit was announced.  We also observe a significant decline between questions 9-11 and a significant rise for the last 2 questions.  This was due to the fact that Contestant $B$ had been getting more questions correct than expected but hit his Lucky Number when he got his remaining 2 Unsure questions incorrect.  Since the probability of getting 2 Unsure questions incorrect in a row is only $0.25^2 = 0.0625$, Contestant $B$ certainly got lucky.

\section{Case Study 2}
\begin{figure}
    \centering
    \includegraphics[width=0.75\linewidth]{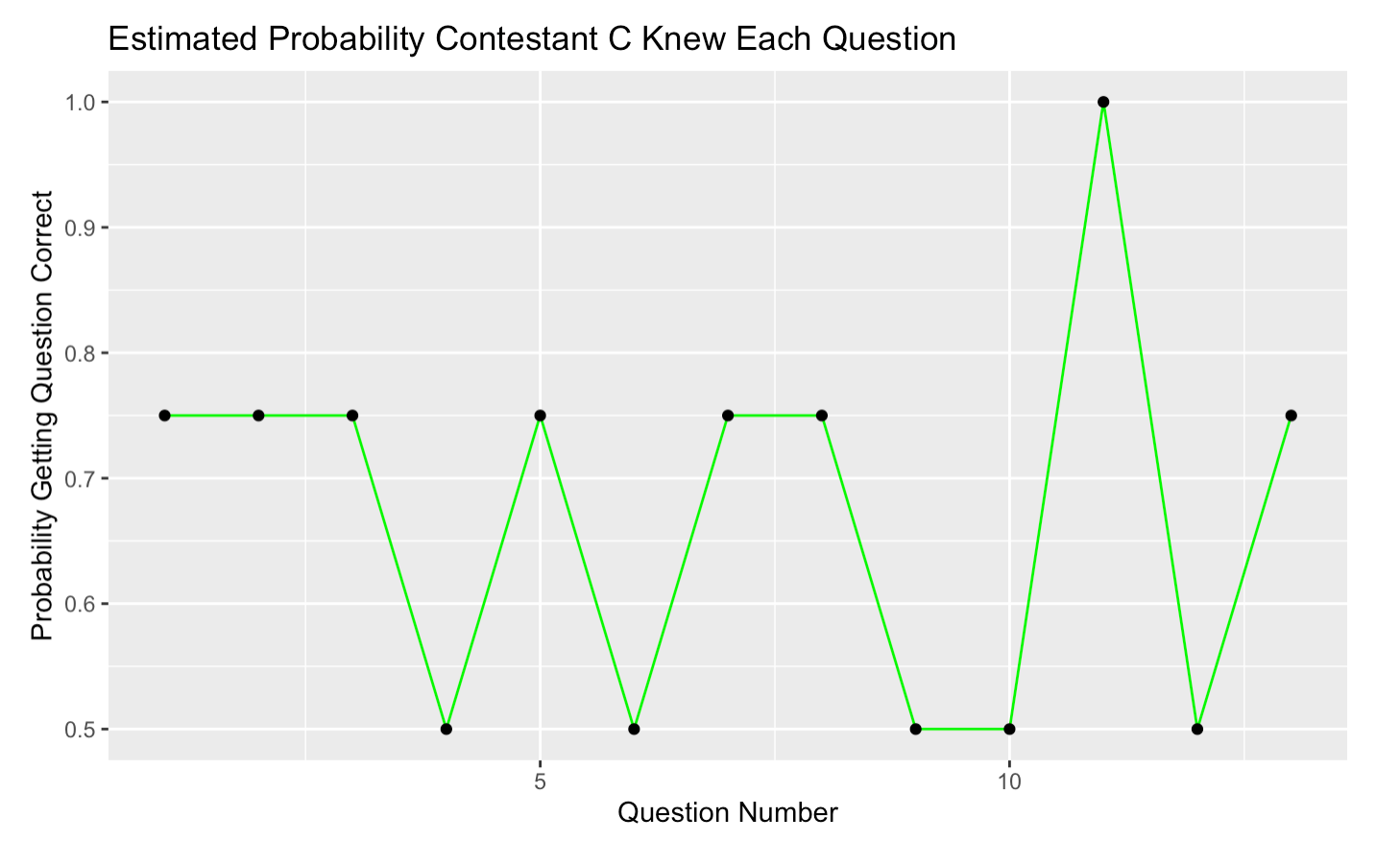}
    \caption{Estimated Probabilities of Correct Answers for Contestant $C$}
    \label{Fig4}
\end{figure}

\begin{figure}
    \centering
    \includegraphics[width=1\linewidth]{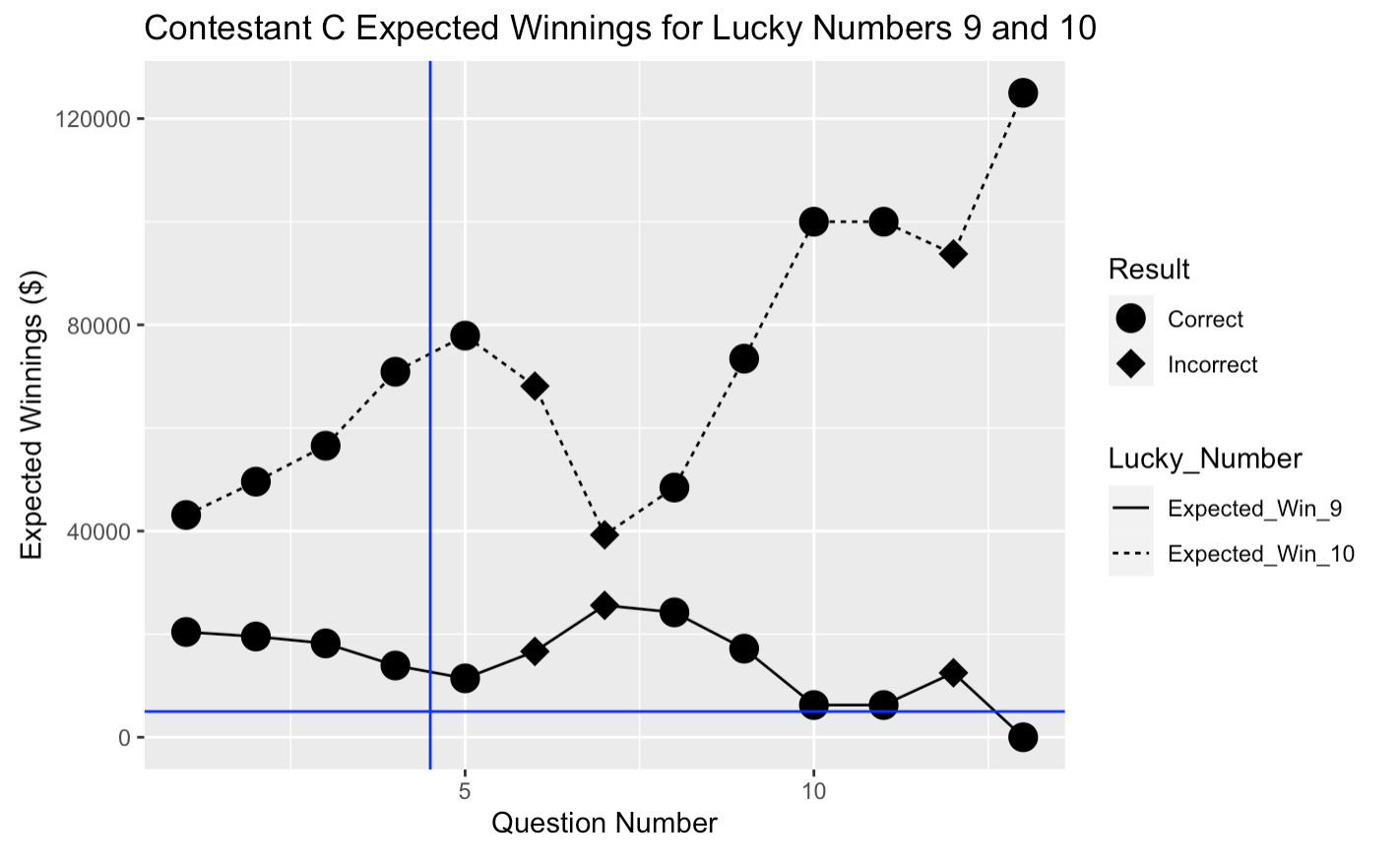}
    \caption{Expected Winnings for Contestant $C$}
    \label{Fig5}
\end{figure}

Now we will examine the game of a second player, Contestant $C$, who selected Lucky Range 7-9 and Lucky Number 9. We determined that she has $S=1, U=7$, and $G=5$. Figure \ref{Fig4} shows the order in which the producers chose to reveal her answers.  Contestant $C$'s expected winnings before any answers were revealed are calculated to be approximately \$20,866.  Figure \ref{Fig5} shows the trajectory of Contestant $C$'s expected winnings after each correct answer is revealed.  After the answer to question 4 was revealed, Contestant $C$ had 4 correct answers, $S=1$, $U=4$, and $G=4$.  At this point, she rejected the offer of \$5000 from the producers.  Under our model, her expected winnings at this point were approximately \$13,916 so it was the correct decision to reject the offer at the time.  The probability Contestant $C$ is in her Lucky Range is 0.34 and the probability she hits her Lucky Number is 0.22.  Unfortunately, she ended up getting 10 questions correct and won nothing.  We argue that contestant $C$ would have been better off to select Lucky Range 10-12 with Lucky Number 10 where she would have won \$125,000.  Figure \ref{Fig5} also shows the trajectory of expected winnings with this new Lucky Range.  With the selection of Lucky Number 10, her starting expectation of winnings increases from approximately \$20,866 to approximately \$37,210.  In this game, knowing your confidence level on questions is more important than knowing the correct answers.

\section{Simulation of Number of Correct Responses}
Next, we develop a model to determine the number of questions an average player is likely to get correct.  After watching the games of 12 contestants (156 questions), we separate the questions into categories (one category per question) based on subjective judgment.  We use this as data to estimate the true proportion of questions falling into each category.  The results are shown in Table \ref{Tab7}.

\begin{figure}
    \centering
    \includegraphics[width=0.8\linewidth]{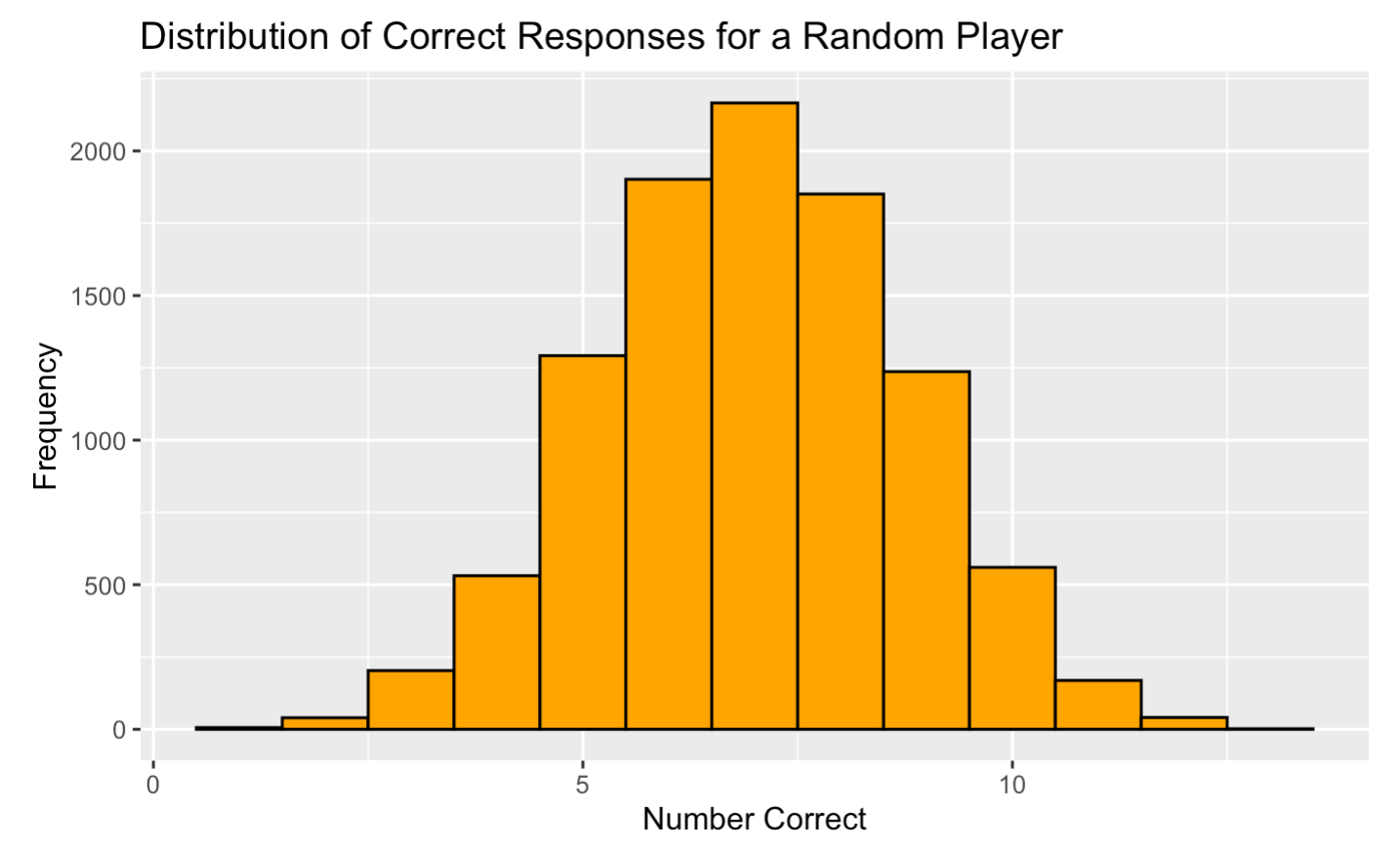}
    \caption{Distribution of Correct Responses}
    \label{Fig6}
\end{figure}

\begin{table}[h]
\caption{Table of Question Categories}\label{Tab7}%
\begin{tabular}{@{}ll@{}}
\toprule
Category & Est. Prob. of Seeing This Category \\
\midrule
Celebrities & 0.12 \\
Movies/TV & 0.08 \\
Animals/Biology & 0.08 \\
History & 0.07 \\
Sports & 0.07 \\
Geography & 0.07 \\
Words & 0.05 \\
Musicians & 0.05 \\
Food & 0.04 \\
U.S. Politicians & 0.04 \\
Inventions & 0.03 \\
U.S. States & 0.03 \\
Space & 0.03 \\
U.S. Gov./Laws & 0.03 \\
Holidays & 0.02 \\
Literature/Magazines & 0.02 \\
Theater & 0.02 \\
Landmarks & 0.02 \\
Periodic Table of Elements & 0.02 \\
Business & 0.02 \\
Other & 0.10 \\
\botrule
\end{tabular}
\end{table}
We assume that a player that is randomly selected has an expertise in a random number of these categories.  If the player has an expertise in the question category, we assume that question will be $U$ with probability 0.6 and $S$ with probability 0.4 which would make the total probability the player gives a correct answer 0.85.  Otherwise, we assume the player will supply the correct answer with probability 0.5.  The number of categories a randomly selected player has an expertise on is modeled by $Q:= R + 1$ where $R \sim \text{Binom}(20, 1.5/20)$. The parameter of 20 was selected since there were 20 categories (excluding the ``Other" category) of questions in Table \ref{Tab7} and the success probability was chosen so that on average a player has expertise in 2.5 categories. We sample without replacement from the set of question categories given in Table \ref{Tab7} in order to select which categories a player will have an expertise in.  For each of the 10,000 players in our Monte Carlo simulation we compute $Q$, take a sample of size $Q$ from the set of all categories, and have the player answer 13 randomly selected questions with the modified probabilities discussed above if a player has an expertise in a category.  The distribution of correct responses under this model is given in Figure \ref{Fig6}.  We observe that the most frequently occurring numbers of correct responses are 6,7, or 8.  Thus, it would be wise for \textit{Lucky 13} contestants to memorize the correct Lucky Number based on their utility function for values $S,U,$ and $G$ such that $6 \leq \text{E}[N] \leq 8$. We also note that it is possible to determine the exact distribution of the number of correct answers but we chose to do a Monte Carlo simulation instead.  

\section{Extended Three Category Model}
The major flaw with the previously discussed three-category model is that contestants may have significant variability in the estimation of the probability that they answer an Unsure question correctly.  In order to study this scenario, we will leave the Sure and Guess categories unchanged and modify the model for the Unsure category. Now, for all questions that are not Sure or Guess, we will assign each of them probability $\pi$ of being a Sure question and probability $1-\pi$ of being a Guess question. We can then impose a prior probability distribution on $\pi$ which leads us to the following theorem.

\begin{theorem}
    If $X_i \sim \text{Bern}(\pi_i)$ and $\pi_i \stackrel{iid}{\sim} F$ with CDF $F(\pi)$ being any probability distribution on $[0,1]$ with PDF $f(\pi)$, then the marginal distribution of $X_i$ is $\text{Bern}(\text{E}[\pi_i])$.
\end{theorem}

\begin{proof}
  Let $N_U$ be the number of questions answered correctly out of $U$ Unsure questions. Clearly, $N_U = \sum_{i=1}^UX_i$ where $X_i$ is the outcome of the $i$-th question. Let $\pi_i$ be an independent sample from $f(\pi)$.  Thus,
  $$\text{P}(X_i=1) = \int_0^1 \text{P}(X_i=1|\pi_i)f(\pi_i)d\pi_i$$
  $$=\int_0^1 \pi_if(\pi_i)d\pi_i$$
  $$= \text{E}[\pi].$$
  Since $\pi_i$ are independent, $X_i$ are marginally independent Bernoulli trials with success probability $\text{E}[\pi]$.  Thus, $N_U \sim \text{Binom}(U,\text{E}[\pi])$.  
\end{proof}

\begin{coro}
Any two prior distributions with the same mean will yield the same strategy.
\end{coro}

\begin{proof}
Since the values of $N_S$ and $N_G$ did not depend on $\pi_i$, any two prior distributions with the same mean will yield the same values for $\text{P}(N=k|\vec{\pi})$ where $k \in [0,13]$ and $\vec{\pi}= (\pi_1,\cdots, \pi_U)$.  
\end{proof}

\begin{figure}[h!]
    \centering
    \includegraphics[width=0.8\linewidth]{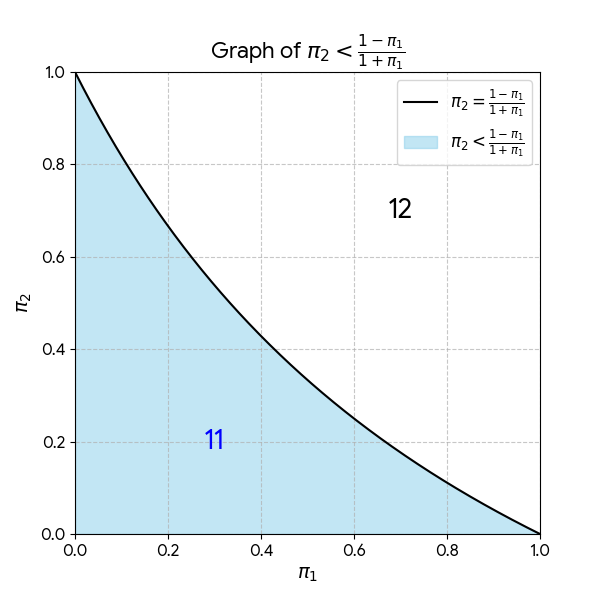}
    \caption{Impact of $\vec{\pi}$ when $S=9, U=2, \text{ and } G=2$}
    \label{Fig7}
\end{figure}

The uniform distribution is a natural choice for the prior distribution of $\pi_i$.  For this choice of prior, we can easily calculate the expected value and variance of 
$$N = S + \text{Binom}\left(G, \frac{1}{2}\right)+\sum_{i=1}^U\text{Bern}\left( \pi_i\cdot 1 + (1-\pi_i)\cdot \frac{1}{2}\right)$$
$$= S + \text{Binom}\left(G, \frac{1}{2}\right)+\sum_{i=1}^U\text{Bern}\left( \frac{\pi_i+1}{2}\right).$$
Assuming $\pi_i \sim \text{Unif}[a_i,b_i]$ for real numbers $a_i, b_i$ such that $0\leq a_i \leq b_i \leq 1$,  we find that
$$\text{E}[N] = \text{E}[\text{E}[N|\vec{\pi}]]$$
$$ = \text{E}\left[S + \frac{G}{2}+ \sum_{i=1}^U\text{E}\left[\text{Bern}\left( \frac{\pi_i+1}{2}\right)\right]\right]$$
$$ = S + \frac{G}{2}+ \sum_{i=1}^U \text{E}\left[\frac{\pi_i+1}{2}\right]$$
$$ = S + \frac{G}{2}+\sum_{i=1}^U\left(\frac{1}{2}\cdot \frac{a_i+b_i}{2}+\frac{1}{2} \right)$$
$$ = S + \frac{G}{2}+\frac{U}{2}+\frac{1}{4}\sum_{i=1}^U(a_i+b_i)$$
and
$$\text{Var}[N] = \text{E}[\text{Var}[N|\vec{\pi}]]+\text{Var}[\text{E}[N|\vec{\pi}]]$$
$$ = \text{E}\left[0+\frac{G}{4}+\sum_{i=1}^U\text{Var}\left[\text{Bern}\left(\frac{\pi_i+1}{2} \right) \right] \right]+\text{Var}\left[S + \frac{G}{2}+\sum_{i=1}^U\left( \frac{\pi_i+1}{2}\right) \right]$$
$$=\frac{G}{4}+\sum_{i=1}^U\text{E}\left[\left(\frac{\pi_i}{2}+\frac{1}{2}\right)\left(\frac{1}{2}-\frac{\pi_i}{2} \right) \right] + \frac{1}{4}\sum_{i=1}^U\text{Var}[\pi_i]$$
$$ = \frac{G}{4}+\frac{1}{4}\sum_{i=1}^U\text{E}[1-\pi_i^2]+ \frac{1}{48}\sum_{i=1}^U(b_i-a_i)^2$$
$$= \frac{G}{4}+\frac{U}{4}-\frac{1}{4}\sum_{i=1}^U\left(\frac{(b_i-a_i)^2}{12}+\frac{(a_i+b_i)^2}{4} \right)+ \frac{1}{48}\sum_{i=1}^U(b_i-a_i)^2$$
$$= \frac{G}{4}+\frac{U}{4}-\frac{1}{16}\sum_{i=1}^U(a_i+b_i)^2.$$
We now demonstrate the impact of $\vec{\pi}$ on a Lucky Number decision.  Assume a player has $S=9, U=2, \text{ and } G=2$.  Assuming we are also provided with the vector $(\pi_1,\pi_2)$, we can compute the following probabilities:
$$\text{P}(N=13|\vec{\pi}) = \text{P}(N_U=2|\vec{\pi})\text{P}(N_G=2)$$
$$ = \left(\frac{\pi_1+1}{2} \right)\left(\frac{\pi_2+1}{2} \right){2 \choose 0}\left( \frac{1}{2}\right)^2$$
$$= \frac{1}{16}[\pi_1\pi_2+\pi_1+\pi_2+1]$$
\vspace{2mm}
$$\text{P}(N=12|\vec{\pi}) = \text{P}(N_U=2|\vec{\pi})\text{P}(N_G=1)+\text{P}(N_U=1|\vec{\pi})\text{P}(N_G=2)$$
$$= \left(\frac{\pi_1+1}{2} \right)\left(\frac{\pi_2+1}{2} \right){2 \choose 1}\left(\frac{1}{2}\right)^2 + \left(\frac{\pi_1+1}{2} \right)\left(\frac{1-\pi_2}{2} \right){2 \choose 2}\left(\frac{1}{2}\right)^2+ \left(\frac{1-\pi_1}{2} \right)\left(\frac{\pi_2+1}{2} \right)\left(\frac{1}{2}\right)^2$$
$$ = \frac{1}{16}[2\pi_1+2\pi_2+4]$$
\vspace{2mm}
$$\text{P}(N=11|\vec{\pi}) = \text{P}(N_U=2|\vec{\pi})\text{P}(N_G=0)+\text{P}(N_U=0|\vec{\pi})\text{P}(N_G=2)+\text{P}(N_U=1|\vec{\pi})\text{P}(N_G=1)$$
$$ = \frac{1}{16}[6-2\pi_1\pi_2]$$
\vspace{2mm}
$$\text{P}(N=10|\vec{\pi})=\text{P}(N_U=1|\vec{\pi})\text{P}(N_G=0)+\text{P}(N_U=0|\vec{\pi})\text{P}(N_G=1)$$
$$=\frac{1}{16}[4-2\pi_2-2\pi_1]$$
\vspace{2mm}
$$\text{P}(N=9|\vec{\pi})=\text{P}(N_U=0|\vec{\pi})\text{P}(N_G=0)$$
$$=\frac{1}{16}[1-\pi_1-\pi_2+\pi_1\pi_2].$$
Now, we can investigate which Lucky Range and Lucky Number this player should choose in order to maximize the probability of winning.  In order for 13 to have maximum probability over Lucky Range 10-12, we would need to have:
$$\text{P}(N=13|\vec{\pi}) > \sum_{i=10}^{12}\text{P}(N=i|\vec{\pi})$$
which is equivalent to 
$$\pi_1+\pi_2+\pi_1\pi_2+1 > 2\pi_1+2\pi_2+4+6-2\pi_1\pi_2+4-2\pi_2-2\pi_1$$
or
$$\pi_1+\pi_2+3\pi_1\pi_2 > 13$$
which clearly cannot happen since $\pi_1,\pi_2 \in [0,1]$.  Thus, claiming all questions were answered correctly is not a better strategy than selecting Lucky Range 10-12.  Now we compare Lucky Ranges 10-12 and 7-9.  In order for Lucky Range 10-12 to have the maximum probability, we would need 
$$\sum_{i=10}^{12}\text{P}(N=i|\vec{\pi}) > \text{P}(N=9|\vec{\pi})$$
which is equivalent to 
$$2\pi_1+2\pi_2+4+6-2\pi_1\pi_2+4-2\pi_2-2\pi_1 > 1-\pi_2-\pi_1+\pi_1\pi_2$$
or
$$13 > 3\pi_1\pi_2-\pi_2-\pi_1.$$
This is always true.  In other words, to maximize the probability of winning, a contestant with $S=9, U=2, \text{ and } G=2$ should always select Lucky Range 10-12.  Now we compare the probabilities of each Lucky Number in the range 10-12.  Note that $P_{10}:= \text{P}(N=10|\vec{\pi})-4=-2\pi_2-2\pi_1 $ and $P_{11}:= \text{P}(N=11|\vec{\pi})-4=2(1-\pi_1\pi_2) $ and $P_{12}:= \text{P}(N=12|\vec{\pi})-4=2(\pi_1+\pi_2) $.  It is clear that $P_{10} < \min\{P_{11},P_{12}\}$ so we only need to compare $P_{11}$ and $P_{12}$.  If $1-\pi_1\pi_2 > \pi_1+\pi_2$ we will select Lucky Number 11.  However, this condition is equivalent to $(1-\pi_1)-\pi_2(1+\pi_1) > 0$ which implies $\pi_2 < \frac{1-\pi_1}{1+\pi_1}$.  Figure \ref{Fig7} shows the regions where we pick Lucky Number 11 versus Lucky Number 12.  

The area of the region shaded in Figure \ref{Fig7} is given by 
$$\int_0^1\frac{1-\pi_1}{1+\pi_1}d\pi_1 = 2\ln(2)-1\approx 0.38.$$
This probability represents the proportion of times that a player would select 11 assuming that the prior distribution of $\pi_1$ and $\pi_2$ are Unif[0,1].
Additionally, Table \ref{Tab8} shows the corresponding values of $S$, $U$, and $G$, for various points in Figure \ref{Fig7}.

\begin{table}[h!]
\caption{Table of Coordinate Point Equivalencies for $S=9, U=2, \text{ and } G=2$}\label{Tab8}%
\begin{tabular}{@{}ll@{}}
\toprule
Point & $S/U/G$ \\
\midrule
$(0,0)$ & 9/0/4\\
$(1,0)$ & 10/0/3\\
$(0,1)$ & 10/0/3\\
$(1,1)$ & 11/0/2\\

\botrule
\end{tabular}
\end{table}

\section{Model with Zero Question Categories}
The most general model for this game is one where the contestant may specify $\vec{p} = (p_1,p_2,\cdots, p_{13})$ where $p_i$ is the probability you get the ith question correct.  If we again define $N$ as the total number of correct responses a contestant gives, $N$ has a Poisson binomial distribution with expected value $\mu = \sum_{i=1}^{13}p_i$.  See \citep{Bib5} for more details on the Poisson binomial distribution.  Note that instead of categories of questions, this model assigns any probability $0.5 \leq p_i \leq 1$ to each question.  Tang and Tang (see \citep{Bib5}) give Darroch's rule for the mode ($m$) of a Poisson binomial random variable based on its mean $\mu$:
\[ m = \begin{cases} 
          k & \text{ if }k \leq \mu < k + \frac{1}{k+2} \\
          k \text{ or } k+1 & \text{ if }k+\frac{1}{k+2} \leq \mu \leq k+1-\frac{1}{n-k+1}\\
          k+1 &\text{ if } k+1-\frac{1}{n-k+1} < \mu \leq k+1.
       \end{cases}
    \]
For any Poisson binomial random variable, $m$ differs from $\mu$ by at most 1.  Darroch's rule allows us to find ranges for $\mu$ that allow us to select the correct Lucky Number in order to maximize the probability of winning as shown in Table \ref{Tab8}.  Further questions about how to select the Lucky Range or when to accept/reject the show's offer for the zero question model would depend on the specific vector $\vec{p}$ and not only $\mu$.  We do not address these questions in this paper but they will be investigated in future research.

\begin{table}[h!]
\caption{Table of Modes}\label{Tab8}%
\begin{tabular}{@{}ll@{}}
\toprule
Mode & Range of $\mu$ \\
\midrule
13 & (12.5,13] \\
12 or 13 & [12.071,12.5]\\
12 & (11.667,12.071)\\
11 or 12 & [11.077,11.667]\\
11 & (10.75,11.077)\\
10 or 11 & [10.0833,10.75]\\
10 & (9.8,10.0833)\\
9 or 10 & [9.091,9.8]\\
9 & (8.833,9.091)\\
8 or 9 & [8.1,8.833]\\
8 & (7.857,8.1)\\
7 or 8 & [7.111,7.857]\\
7 & (6.875,7.111)\\
6 or 7 & [6.125,6.875]\\
6 & (5.888,6.125)\\
\botrule
\end{tabular}
\end{table}

\section{Conclusion and Future Work}
We have presented models for the optimal strategy of \textit{Lucky 13} with 2 and 3 question categories based on two different utility functions.  We have also presented two case studies of actual games and used Monte Carlo simulation to estimate the number of questions an average player may answer correctly.  Additionally, we explain how this model could be extended to allow contestants to provide individual probabilities of giving a correct answer for each of the Unsure questions or all 13 questions.  This research could be extended by conducting an empirical study to determine if contestant decisions in a real game came close to our optimal strategy while playing.  Furthermore, a Monte Carlo simulation could be conducted to estimate the probabilities (and their standard errors) for each Lucky Number under the extended three category model given the contestant data and an appropriate prior distribution.  Both of these ideas are subject to ongoing research.  It would also be interesting to model the amount of money the producers offer a contestant to quit.  However, the amount of data available was too sparse to come up with a realistic model.  \\
\\

\section{Statements and Declarations}

\subsection{Author Contribution}
The manuscript was prepared by S.B. and all authors reviewed the manuscript and results.  
\subsection{Funding}
The authors have no relevant financial or non-financial interests to disclose. 

\subsection{Ethics Approval}
Not applicable.  This research does not qualify as human subject research.

\pagebreak

\begin{appendices}

\section{Optimal Strategy for Selecting Lucky Range and Number (3 Category Model)}\label{secA1}
The first chart shows how to select the correct Lucky Range and Lucky Number if you are playing to maximize the probability you win some money.  The second chart shows how to make your decision if you are playing to maximize your winnings.  Both charts are colored by Lucky Range.

\begin{figure}[h!]
    \centering
    \includegraphics[width=0.92\linewidth]{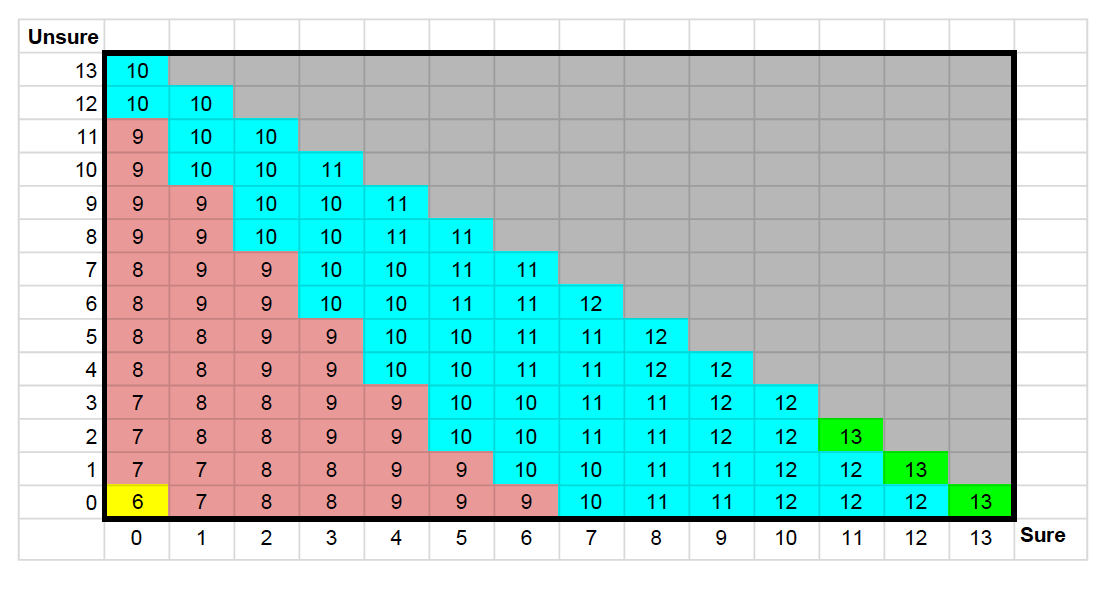}
    \caption{Maximizing Win Probability }
\label{Fig8}
\end{figure}

\begin{figure}[h!]
    \centering

\includegraphics[width=0.92\linewidth]{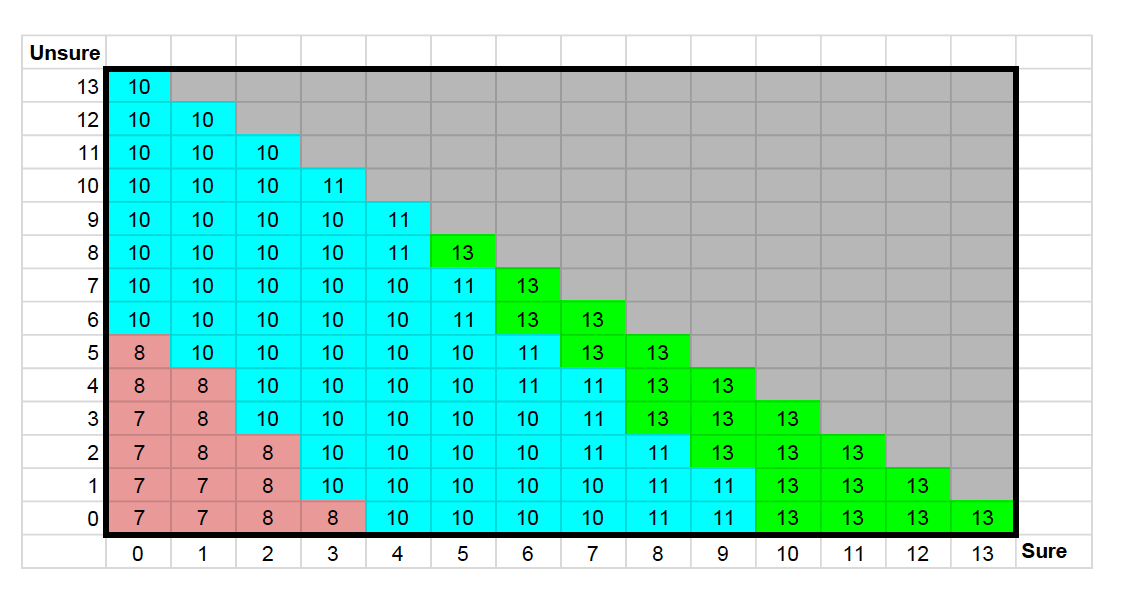}
\caption{Maximizing Expected Winnings}
\label{Fig9}
\end{figure}





\end{appendices}

\pagebreak 
\bibliography{sn-bibliography}

\end{document}